\documentclass{amsart}
\usepackage{amssymb,latexsym}
\theoremstyle{plain}
\newtheorem{theorem}{Theorem}
\newtheorem{corollary}{Corollary}

\newtheorem{lemma}{Lemma}
\theoremstyle{definition}

\newtheorem{remark}{Remark}

\date{}

\begin{document}

\title[regularity]
{$n$-blocks collections on Fano manifolds and sheaves with regularity
 $-\infty$}
\author{E. Ballico}
\address{Dept. of Mathematics\\
 University of Trento\\
38050 Povo (TN), Italy}
\email{ballico@science.unitn.it}
\thanks{This author was partially supported by MIUR and GNSAGA of INdAM
 (Italy).}
\author{F. Malaspina}
\address{Universit\`{a} di Torino\\
via C. Alberto 10, 10123 Torino, Italy}
\email{francesco.malaspina@unito.it}
\subjclass{14J60}
\keywords{$m$-block collection; Fano manifold; coherent sheaves on a
 Fano manifold; order of regularity}

\begin{abstract}
Let $X$ be a smooth Fano manifold equipped with a `` nice '' $n$-blocks
 collection in the sense
of \cite{cm2} and $\mathcal {F}$ a coherent sheaf on $X$. Assume that
 $X$ is Fano and that all blocks are
coherent sheaves. Here we prove that $\mathcal {F}$ has regularity
 $-\infty$ in the sense of \cite{cm2}
if $\mbox{Supp}(\mathcal {F})$ is finite, the converse being true under
 mild assumptions. The corresponding result is also
true when $X$ has a geometric collection in the sense of \cite{cm1}.
\end{abstract}

\maketitle

\section{Introduction}\label{S1}

Let $X$ be an $n$-dimensional smooth projective variety over $\mathbb
 {C}$. Let $\mathcal {D}:= \mathcal {D}^b(\mathcal {O}_X-{\it mod})$
denote the bounded category of $\mathcal {O}_X$-sheaves. Let $\mathcal
 {F}$ be a coherent sheaf
on $X$. Assume that $X$ has a geometric collection
in the sense of \cite{cm1} or an $n$-blocks collection in the sense of
 \cite{cm2}. L. Costa and R. M. Mir\'{o}-Roig
defined the notion of regularity for $\mathcal {F}$ and asked a
 characterization of all $\mathcal {F}$ whose regularity
is $-\infty$ (\cite{cm1}, Remark 3.3). In section \ref{S2} we will
 recall the definitions contained in \cite{cm1} and \cite{cm2}
and used in our statements below. After the statements we will discuss
 our motivations and give a very short list of interesting
varieties to which these results may be applied.

We prove the following results.

\begin{theorem}\label{i1}
Assume that $X$ is Fano and that it has an $n$-blocks collection
 $\mathcal {B}$ whose members are coherent sheaves.
Let $\mathcal {F}$ be a coherent sheaf on $X$. If $\mathcal {F}$ has
 regularity $-\infty$ with respect to $\mathcal {B}$,
then $\mbox{Supp}(\mathcal {F})$ is finite. If all right mutations of
 all elements of $\mathcal {B}$ are locally free
and $\mbox{Supp}(\mathcal {F})$ is finite, then  $\mathcal {F}$ has
 regularity $-\infty$ with respect to $\mathcal {B}$. 
\end{theorem}

\begin{corollary}\label{i2}
Assume that $X$ is Fano and that it has a geometric collection
 $\mathcal {G}$ whose members are coherent sheaves.
Let $\mathcal {F}$ be a coherent sheaf on $X$. If $\mathcal {F}$ has
 regularity $-\infty$ with respect to $\mathcal {G}$,
then $\mbox{Supp}(\mathcal {F})$ is finite. If all right mutations of
 all elements of $\mathcal {G}$ are locally free
and $\mbox{Supp}(\mathcal {F})$ is finite, then  $\mathcal {F}$ has
 regularity $-\infty$ with respect to $\mathcal {G}$. 
\end{corollary}

We recall that any projective manifold with a geometric collection is
 Fano (\cite{cm1}, part (2) of Remark 2.16).
Any $n$-dimensional smooth quadric $Q_n \subset
{\bf {P}}^{n+1}$ has an $n$-block collection whose members are locally
 free (\cite{cm2}, Example 3.2 (2)). It has a geometric collection
if and only if $n$ is odd. Any Grassmannian $G$ has an $n$-block
 collection (with $n:= \dim (G)$) whose members
are locally free sheaves (\cite{cm2}, Example 3.7 (4)). For the Fano
 $3$-folds $V_5$ and $V_{22}$ D. Faenzi found a geometric
collection whose members are locally free
(\cite{f1}, \cite{f2}).

Castelnuovo-Mumford regularity was introduced by Mumford in \cite{m},
 Lecture 14, for a coherent sheaf $\mathcal {F}$ on ${\bf {P}}^n$.
He ascribed the idea to Castelnuovo for the following reason. Let $C
 \subset {\bf {P}}^n$ be a closed subvariety
and $H \subset {\bf {P}}^n$ be a general hyperplane. Then we have an
 exact sequence
\begin{equation}\label{eqi1}
0 \to \mathcal {I}_C(t-1) \to \mathcal {I}_C(t) \to \mathcal {I}_{C\cap
 H}(t) \to 0
\end{equation}
Castelnuovo used the corresponding classical (pre-sheaves) concepts of
 linear systems to get informations on $C$ from informations on $C\cap
 H$
plus other geometrical or numerical assumptions on $C$. The key
 properties of Castelnuovo-Mumford regularity is that if $\mathcal {F}$ is
 $m$-regular,
then it is $(m+1)$-regular and $\mathcal {F}(m)$ (or $\mathcal
 {I}_C(m)$) is spanned. Since \cite{m} several hundred
papers
studied this notion, which is now also a key property in computational
 algebra. 
Let $X$ be a projective scheme, $H$ an ample line bundle on $X$ and
 $\mathcal {F}$ a coherent sheaf on $X$. The definition
in \cite{m}, Lecture 14, apply verbatim, just writing $\mathcal
 {F}\otimes H^{\otimes t}$ instead of $\mathcal {F}(t)$.
This is also called Castelnuovo-Mumford regularity with respect to the
 polarized pair $(X,H)$. $X$ may have several non-proportional
polarizations. It is better to collect all informations for all
 polarizations in a single integer (the regularity)
not in a string of integers, one for each proportional class of
 polarizations on $X$.
This is the reason for the definitions given by Hoffman-Wang for
 products of projective varieties (\cite{hw})
and by Maclagan and Smith for toric varieties (\cite{ms}). Even when
 $X$ has only one polarization
the search for generalizations of Beilinson's spectral sequence from
 ${\bf {P}}^n$ to $X$ gave a strong motivation
to introduce the notions of regularities for geometric collections
 (\cite{cm1}, Th. 2.21) and $n$-block collections
(\cite{cm2}, Th. 3.10).
The reader will notice that to prove Theorem \ref{i1}
and Corollary \ref{i2} we will use neither the main definitions of
 \cite{cm1} and \cite{cm2} nor the machinery of derived categories. We will
 only use the formal properties (like `` spannedness '' or ``
 $m$-regularity implies
$(m+1)$-regularity '') proved in \cite{cm1} and \cite{cm2}
(see eq. (\ref{eqi0}) in section \S 2 for an explanation of the word ``
 spannedness ''). We hope that our results will be
extended and used if
other notions of regularity will appear in the literature.

\section{The main definitions and the proofs}\label{S2}

Let $X$ be an $n$-dimensional smooth projective variety over $\mathbb
 {C}$. Let $\mathcal {D}:= \mathcal {D}^b(\mathcal {O}_X-{\it mod})$
denote the bounded category of $\mathcal {O}_X$-sheaves. For all
 objects $A, B\in \mathcal {D}$
set $\mbox{Hom}^{\bullet }(A,B):= \oplus _{k\in \mathbb {Z}}
 \mbox{Ext}^k_\mathcal {D}(A,B)$.
An object $A\in \mathcal {D}$
is said to be {\it exceptional} if $\mbox{Hom}^{\bullet }(A,A)$ is an
 $1$-dimensional algebra generated by the identity. An ordered
collection $(A_0,\dots ,A_m)$ of objects of $\mathcal {D}$ will be
 called an {\it exceptional collection} if each
$A_i$ is exceptional and $\mbox{Ext}^{\bullet}_\mathcal {D}(A_k,A_j) =
 0$ for all $0 \le j < k \le m$. A collection
$(A_0,\dots ,A_m)$ is said to be {\it strongly exceptional} if it is
 exceptional and $\mbox{Ext}^i_\mathcal {D}(A_j,A_k)=0$
for all $(i,j,k)$ such that $i \ne 0$ and $j \le k$. A collection
 $(A_0,\dots ,A_m)$ is said to be full
if it generates $\mathcal {D}$. This implies $\mathcal {D} \cong
 \mathbb {Z}^{\oplus (m+1)}$. Now asssume that $X$
admits a fully exceptional collection $\sigma = (A_0,\dots ,A_n)$. For
 any $A, B\in \mathcal {D}$
the right mutation $R_BA$ of $A$ and the left mutation $L_AB$ of $B$
 are defined by the following distinguished triangles
\begin{equation*}
R_BA[-1] \to A \to \mbox{Hom}^{\times \bullet}(A,B)\otimes B \to R_BA
\end{equation*}
\begin{equation*}
L_AB \to \mbox{Hom}^\bullet (A,B)\otimes A \to L_AB[1]
\end{equation*}
(\cite{cm1}, Definition 2.4). For every integer $i$ such that $1 \le i
 \le n$, define the $i$-th right mutation $R_i\sigma$
and the $i$-th left mutation $L_i\sigma$ of $\sigma$
by the formulas
$$R_i\sigma := (A_0,\dots A_{i-2},A_i,R_{A_i}A_{i-1},A_{i+1},\dots
 A_n)$$
$$L_i\sigma := (A_0,\dots ,A_{i-2},
 L_{A_{i-1}}A_i,A_{i-1},A_{i+1},\dots ,A_n)$$
(a switch of two elements of $\sigma$ and the application
to one of them of a right or left mutation) 
(\cite{cm1}, Definition 2.6). For any $j \ge 2$, set $R^{(j)}A_i:=
 R_{A_{i+j}}\circ \cdots \circ R_{A_{i+1}}A_i\in \mathcal {D}$ and
 define in a similar way the
iterated left mutations $L^{(i)}$ (\cite{cm1}, Notation 2.7). Set
 $A_{n+i}:= R^{(n)}A_{i-1}$
for all $0 \le i \le n$
and $A_{-i}:= L^{(n)}A_{n-i+1}$ for all $1 \le i \le n$. Iterating the
 use of $R^{(n)}$ and $L^{(n)}$ we
get the elix $\{A_i\}_{i\in \mathbb {Z}}$ with $A_i\in \mathcal {D}$
 for all $i$
(\cite{cm1}, Definition 2.12). For instance, if $X={\bf {P}}^n$, then
$(A_0,\dots ,A_n):= (\mathcal {O}_{{\bf {P}}^n},\mathcal {O}_{{\bf
 {P}}^n}(1),\dots ,\mathcal {O}_{{\bf {P}}^n}(n))$
is a geometric collection and $\{\mathcal {O}_{{\bf {P}}^n}(t)\}_{t\in
 \mathbb {Z}}$ is the corresponding elix.
Let $\mathcal {F}$ be a coherent
sheaf on $X$. $\mathcal {F}$ is said to be $m$-regular with respect to
 the geometric collection $\sigma = (A_0,\dots ,A_n)$
if $\mbox{Ext}^q(R^{(-p)}A_{-m+p},\mathcal {F}) =0$ for all integers
 $q,p$ such that $q>0$ and $-n \le p \le 0$. The regularity of $\mathcal
 {F}$ is the minimal integer $m$ such that
$\mathcal {F}$ is $m$-regular (or $-\infty$ if it is $m$-regular for
 all $m\in \mathbb {Z}$).
An exceptional collection $(A_0,\dots ,A_s)$ of elements of $\mathcal
 {D}$ is called a {\it block} if $\mbox{Ext}^i_\mathcal {D}(
A_j,A_k)=0$ for all $i, j, k$ such that $k\ne j$. An $m$-block
 collection of elements of $\mathcal {D}$
is an exceptional collection which may be partitioned into $m+1$
 consecutive blocks. Assume that $X$ has an $n$-block collection
whose elements generate $\mathcal {D}$. Let $\mathcal {F}$ be a
 coherent sheaf on $X$. In \cite{cm2}, Definition 4.5,
there is a definition of regularity of $\mathcal {F}$; it requires only
 technical modifications with
respect to the simpler case of a geometric collection: they gave
 similar definitions of left and right mutations and elices.
Then the definition of $m$-regularity is again given by certain
 $\mbox{Ext}$-vanishings. If a coherent sheaf $\mathcal {F}$
is $m$-regular with respect to a geometric collection $\sigma$ or an
 $n$-block collection $\sigma$, then it gives a resolution
\begin{equation}\label{eqi0}
0 \to \mathcal {L}_{-n} \to \cdots \to \mathcal {L}_{-1} \to \mathcal
 {L}_0 \to \mathcal {F} \to 0
\end{equation}
in which each $\mathcal {L}_i\in \mathcal {D}$ is constructed from
 $\mathcal {F}$ and the elements of $\sigma$
taking tensor products (\cite{cm1}, between 3.1 and 3.2 for geometric
 collections, \cite{cm2}, eq. (4.2), for $n$-blocks).
If the elements of $\sigma$ are coherent sheaves (resp. localy free
 coherent sheaves), then each $\mathcal {L}_i$
is a coherent sheaf (resp. a locally free coherent sheaf). In the case
 of Castelnuovo-Mumford regularity
the corresponding result is true. It shows how the Castelnuovo-Mumford
 regularity bounds the degrees
of the syzygies. This is the key reason for its use in computational
 algebra.

The following well-known result answers the corresponding problem for
 Castelnuovo-Mumford regularity.

\begin{lemma}\label{a1}
Let $X$ be a projective scheme, $L$ an ample line bundle on $X$ and
 $\mathcal {F}$ a coherent sheaf
on $X$. The following conditions are equivalent:
\begin{itemize}
\item[(a)] $\mathcal {F}$ is supported by finitely many points of $X$;
\item[(b)] $\mathcal {F}\otimes L^{\otimes t}$ is spanned for all $t
 \ll 0$;
\item[(c)] $h^i(X,\mathcal {F}\otimes L^{\otimes t})=0$ for all $i>0$
 and all $t\in \mathbb {Z}$.
\end{itemize}
\end{lemma}

\begin{proof}
Obviously, (a) implies (b) and (c). Now assume that (b) holds,
but that $\dim (\mbox{Supp}(\mathcal {F}))>0$. Take an integral
 projective
curve $C \subseteq \mbox{Supp}(\mathcal {F})$. Since the restriction of
 a spanned sheaf
is spanned, $\mathcal {F}\vert C$ satisfies (c) with respect to the
 ample line bundle $R:= L\vert C$.
Let $f: D \to C$ be the normalization. Set $M:= f^\ast (R)$. $M$ is an
 ample
line bundle on $D$. Since $D$ is a smooth curve, the coherent sheaf
$f^\ast (\mathcal {F})$ is either a torsion sheaf or
the direct sum of a torsion sheaf $T$ and a vector bundle $E$ with
 positive rank. To prove (a) we must check
that $f^\ast (\mathcal {F})$ is torsion. Assume $E \ne 0$. Since the
 pull-back of a spanned sheaf is
spanned, $E\otimes M^{\otimes t}$ is spanned for all $t\in \mathbb
 {Z}$. Since $\deg (E\otimes R^{\otimes t})
=\deg (E) + t\cdot \mbox{rank}(E)\cdot \deg (M) < 0$ for $t \ll 0$,
 $E\otimes R^{\otimes t}$ is not
spanned for $t \ll 0$, contradiction. 
Let $x \ge 1$ be an integer such that $L^{\otimes x}$ is very
ample. If $\mathcal {F}$ satisfies (c) for the line bundle $L$, then it
 satisfies
the same condition for the line bundle $L':= L^{\otimes x}$. Hence
to check that (c) implies (a) we may assume that $L$ is very ample. Fix
 an integer $t$. Since
$h^i(X,\mathcal {F}\otimes L^{\otimes (t-i-1)}) = 0$ for all $i>0$,
 $\mathcal {F}\otimes L^{\otimes t}$
is spanned (\cite{m}, p. 100). Thus (b) holds and hence (a) holds.
\end{proof}

\vspace{0.3cm}

\qquad {\emph {Proof of Theorem \ref{i1}.}} Fix a coherent sheaf
 $\mathcal {F}$.
Let $\mathcal {E}$ be the helix of blocks generated by $\mathcal {B}$
(\cite{cm2}, Definition 4.1).
All elements of $\mathcal {E}$ are coherent sheaves, not just complexes
 (\cite{cm2}, Corollary 4.4)
and their elements satisfies a periodicity modulo $n+1$: $\mathcal
 {E}_i = \mathcal {E}_{i+n+1}\otimes \omega _X$
(\cite{cm2}, lines between 4.3 and 4.4). First assume that $\mathcal
 {F}$ has regularity $-\infty$ with respect to $\mathcal {B}$,
i.e. that it is $m$-regular with respect to $\mathcal {B}$ for all $m
 \ll 0$. Fix $m\in \mathbb {Z}$.
The $m$-regularity of $\mathcal {F}$ implies that it is a quotient of a
 finite sum $\mathcal {L}_0$
of sheaves of the form
$E_s^{-m}$ appearing in the blocks of $\mathcal {B}$ (\cite{cm2},
 Definition 4.5). Since
$\mathcal {F}$ is $t$-regular for all $t \ll 0$, the periodicity
property of $\mathcal {E}$ shows that for all integers $t \le 0$,
 $\mathcal {F}$ is a quotient of a finite direct
sum of sheaves of the form $\mathcal {L}_0\otimes \omega _X^{\otimes
 (-t)}$. Since $X$ is Fano, $\omega _X^\ast$
is ample. Take $L:= \omega _X^\ast$ and copy the proof that (b) implies
 (a) in Lemma \ref{a1}. We get that
$\mbox{Supp}(\mathcal {F})$ is finite.\\
Now assume that
$\mbox{Supp}(\mathcal {F})$ is finite and that all right mutations of
 elements of $\mathcal {B}$ are locally free.
Let $A$ be any of these mutations. Since $A$ is locally free, the local
 Ext-functors ${\it Ext}^i(A,\mathcal {F})$ vanish
for all $i>0$. Hence the local-to-global spectral sequence for the
 Ext-functors gives $\mbox{Ext}^i(A,\mathcal {F})
\cong H^i(X,{\it Hom}(A,\mathcal {F}))$ for all $i \ge 0$. Since
 $\mbox{Supp}(\mathcal {F})$ is finite, we get
$\mbox{Ext}^q(A,\mathcal {F}) =0$ for all $q>0$. Hence for every
 integer $m$ the sheaf $\mathcal {F}$
satisfies the definition of $m$-regularity given in \cite{cm2},
 Definition 4.5. Since $\mathcal {F}$ is $m$-regular
with respect to $\mathcal {B}$ for all $m$, its regularity is
 $-\infty$.\qed

\vspace{0.3cm}

\qquad {\emph {Proof of Corollary \ref{i2}.}} This result is a
 particular case of Theorem \ref{i1}, because
the definition of regularity for geometric collections given in
 \cite{cm1} agrees with the definition
of regularity for $n$-blocks collections given in \cite{cm2} (see
 \cite{cm2}, Remark 4.7). It
may be proved directly, just quoting \cite{cm1}, Remark 2.14, to get
 the periodicity
property $\mathcal {E}_i = \mathcal {E}_{i+n+1}\otimes \omega _X$ and
 \cite{cm1}, Proposition 3.8, to get the surjection
$\mathcal {L}_0 \to \mathcal {F}$.\qed 

\begin{remark}\label{a2}
In \cite{c} J. V. Chipalkatti defined a notion of regularity for a
 coherent sheaf $\mathcal {F}$
on a Grassmannian. He remarked
that $\mathcal {F}$ have regularity $-\infty$ (according to his
 definition) if and only if its support
is finite (\cite{c}, part 4) of Remark 1.2).
\end{remark}

\begin{remark}\label{a2}
Let $\mathcal {F}$ be a coherent sheaf on ${\bf {P}}^n\times {\bf
 {P}}^m$. Hoffman and Wang introduced
a bigraded definition of regularity (\cite{hw}). The definition of
 ampleness and \cite{hw}, Prop. 2.8, imply
that if $\mathcal {F}$ is $(a,b)$-regular in the sense of Hoffman-Wang
 for all $(a,b)\in \mathbb {Z}\times \mathbb {Z}$,
then $\mbox{Supp}(\mathcal {F})$ is finite. The converse is obvious. As
 remarked in \cite{cm2}, Remark 5.2, Hoffman-Wang
definition and its main properties may be extended verbatim to
 arbitrary multiprojective spaces ${\bf {P}}^{n_1}\times \cdots
\times {\bf {P}}^{n_s}$.
\end{remark}

\providecommand{\bysame}{\leavevmode\hbox to3em{\hrulefill}\thinspace}


\begin{thebibliography}{99}

\bibitem{c} J. V. Chipalkatti, A generalization of Castelnuovo
 regularity to Grassman varieties, Manuscripta Math.
102 (2000), no. 4, 447--464.

\bibitem{cm1} L. Costa and R. M. Mir\'{o}-Roig, Geometric collections
 and Castelnuovo-Mumford regularity, arXiv:math/0609561,
Math. Proc. Cambridge Phil. Soc., to appear.

\bibitem{cm2} L. Costa and R. M. Mir\'{o}-Roig, $m$-blocks collections
 and Castelnuovo-Mumford regularity
in multiprojective spaces, Nagoya Math. J. 186 (2007), 119-155.

\bibitem{f1} D. Faenzi, Bundles over the Fano threefold $V_5$, Comm.
 Algebra 33 (2005), no. 9, 3061--3080.

\bibitem{f2} D. Faenzi, Bundles over Fano threefolds of type $V_{22}$,
 Ann. Mat. Pura Appl. (4) 186 (2007), no. 1,
1--25.

\bibitem{hw} J. W. Hoffman and H. H. Wang, Castelnuovo-Mumford
 regularity in biprojective spaces, Adv. Geom. 4 (2004),
no. 4, 513--536.

\bibitem{ms} D. Maclagan and G. G. Smith, Multigraded
 Castelnuovo-Mumford regularity, J. Reine Angew. Math. 571
(2004), 179--212.

\bibitem{m} D. Mumford, Lectures on curves on an algebraic surface,
 Princeton University Press, Princeton, N.J., 1966.

\end{thebibliography}
\end{document}